\documentclass[preprint,3p,times]{elsarticle}

\usepackage{amssymb}
\usepackage{amsmath}
\usepackage[all,cmtip]{xy}
\usepackage{amsthm}
\usepackage{color}
\usepackage{enumitem}
\usepackage{tikz}
\usepackage{appendix}
\usepackage[colorlinks=true]{hyperref}
\usepackage[hyphenbreaks]{breakurl}
\input diagxy

\theoremstyle{plain}
  \newtheorem{thm}{Theorem}[section]
  \newtheorem{lem}[thm]{Lemma}
  \newtheorem{prop}[thm]{Proposition}
  \newtheorem{cor}[thm]{Corollary}
\theoremstyle{definition}

  \newtheorem{ques}[thm]{Question}

\DeclareMathAlphabet{\mathcal}{OMS}{cmsy}{m}{n}

\makeatletter
\def\ps@pprintTitle{%
 \let\@oddhead\@empty
 \let\@evenhead\@empty
 \def\@oddfoot{\centerline{\thepage}}%
 \let\@evenfoot\@oddfoot}
\makeatother

\newcommand{\lra}{\longrightarrow}

\newcommand{\CA}{\mathcal{A}}
\newcommand{\CB}{\mathcal{B}}

\newcommand{\CD}{\mathcal{D}}

\newcommand{\CG}{\mathcal{G}}

\newcommand{\CM}{\mathcal{M}}
\newcommand{\CN}{\mathcal{N}}

\newcommand{\CT}{\mathcal{T}}

\newcommand{\sC}{\mathsf{C}}

\newcommand{\cl}{\mathsf{cl}}

\newcommand{\App}{\mathbf{App}}

\newcommand{\Met}{\mathbf{Met}}

\newcommand{\Sup}{\mathbf{Sup}}
\renewcommand{\Top}{\mathbf{Top}}

\newcommand{\ProbMet}{\mathbf{ProbMet}}

\newcommand{\op}{{\rm op}}

\renewcommand{\leq}{\leqslant}
\renewcommand{\geq}{\geqslant}

\renewcommand{\phi}{\varphi}

\numberwithin{equation}{section}

\begin{document}

\begin{frontmatter}



\title{On the probabilistic metrizability of approach spaces}


\author[S]{Hongliang Lai}
\ead{hllai@scu.edu.cn}

\author[S]{Lili Shen}
\ead{shenlili@scu.edu.cn}

\author[B]{Junche Yu\corref{cor}}
\ead{cqyjc@icloud.com}

\cortext[cor]{Corresponding author.}
\address[S]{School of Mathematics, Sichuan University, Chengdu 610064, China}
\address[B]{Laboratoire de Math{\' e}matiques Pures et Appliqu{\' e}es, Universit{\' e} du Littoral-C{\^ o}te d'Opale, 62228 Calais, France}

\begin{abstract}
We investigate approach spaces generated by probabilistic metric spaces with respect to a continuous t-norm $*$ on the unit interval $[0,1]$. Let $k^*$ be the supremum of the idempotent elements of $*$ in $[0,1)$. It is shown that if $k^*=1$ (resp. $k^*<1$), then an approach space is probabilistic metrizable with respect to $*$ if and only if it is probabilistic metrizable with respect to the minimum (resp. product) t-norm. 
\end{abstract}

\begin{keyword}
Approach space \sep Probabilistic metric space \sep Probabilistic metrizable approach space \sep Continuous t-norm

\MSC[2020] 54A05 \sep 54E70 \sep 54E35
\end{keyword}

\end{frontmatter}




\section{Introduction}

\emph{Probabilistic metric spaces} \cite{Menger1942,Schweizer1983} are a generalization of metric spaces in which the distance takes values in distance distributions instead of non-negative real numbers. \emph{Approach spaces} \cite{Lowen1989b,Lowen1997,Lowen2015}, introduced by Lowen, are a common extension for topological spaces and metric spaces. 

The motivation of this paper originates from Lawvere's presentation of metric spaces as enriched categories \cite{Lawvere1973}, and the discovery of Manes and Barr that topological spaces can be encoded in terms of ultrafilter convergence \cite{Manes1969,Barr1970}, which give rise to the theory of \emph{monoidal topology} \cite{Clementino2003,Clementino2004,Hofmann2014}. Explicitly:
\begin{itemize}
\item Classical metric spaces (in the sense of Fr{\' e}chet \cite{Frechet1906}) are symmetric, separated and finitary $[0,\infty]^{\op}$-categories, and probabilistic metric spaces are symmetric, separated and finitary $\Delta^+$-categories \cite{Chai2009,Hofmann2013a}, where $[0,\infty]^{\op}$ is the Lawvere quantale \cite{Lawvere1973,Rosenthal1990}, and $\Delta^+$ is the quantale of distance distributions. 
\item Classical topological spaces are ${\bf 2}$-valued topological spaces, and approach spaces are $[0,\infty]^{\op}$-valued topological spaces \cite{Lai2017a}, where ${\bf 2}$ is the two-element Boolean algebra.
\end{itemize}

Since the Lawvere quantale $[0,\infty]^{\op}=([0,\infty]^{\op},+,0)$ is isomorphic to the quantale $[0,1]=([0,1],\times,1)$ of the unit interval equipped with the usual product, the classical metrizability of topological spaces is actually to consider
\begin{itemize}
\item ${\bf 2}$-valued topological spaces generated by symmetric, separated and finitary $[0,1]$-categories. 
\end{itemize}
So, it is natural to consider its counterpart by taking tensor products $\boxtimes$ in the category $\Sup$ of complete lattices and $\sup$-preserving maps \cite{Joyal1984}; that is, 
\begin{itemize}
\item $([0,\infty]^{\op}\boxtimes{\bf 2})$-valued topological spaces generated by symmetric, separated and finitary $([0,\infty]^{\op}\boxtimes[0,1])$-categories.
\end{itemize}
Since $[0,\infty]^{\op}\boxtimes{\bf 2}=[0,\infty]^{\op}$ (see \cite[Proposition I.5.2]{Joyal1984}) and $[0,\infty]^{\op}\boxtimes[0,1]=\Delta^+$ (see \cite[Subsection 3.3]{GutierrezGarcia2017a}), we are precisely asking for 
\begin{itemize}
\item $[0,\infty]^{\op}$-valued topological spaces generated by symmetric, separated and finitary $\Delta^+$-categories;
\end{itemize}
that is,
\begin{itemize}
\item approach spaces generated by probabilistic metric spaces.
\end{itemize}

As a first step towards the probabilistic metrizability of approach spaces, in this paper we focus on the connections between approach spaces generated by probabilistic metric spaces with respect to different \emph{continuous t-norms} \cite{Klement2000,Alsina2006}. Explicitly, the triangle inequality in a probabilistic metric space $(X,\alpha)$ can be expressed as  
\[\alpha(y,z,r)*\alpha(x,y,s)\leq\alpha(x,z,r+s)\]
for all $x,y,z\in X$, $r,s\in[0,\infty]$, where $*$ is a continuous t-norm\footnote{In some literature (e.g., \cite{Schweizer1983}) the t-norm $*$ is only required to be left-continuous.} on the unit interval $[0,1]$. Therefore, the probabilistic metrizability of an approach space relies on the choice of the continuous t-norm $*$. We say that an approach space is \emph{probabilistic metrizable with respect to $*$}, or \emph{$*$-metrizable} for short, if it is generated by a probabilistic metric with respect to a continuous t-norm $*$ on $[0,1]$. The main result of this paper, Theorem \ref{main}, reveals that for each continuous t-norm $*$ on $[0,1]$, the $*$-metrizability of approach spaces depends on the supremum of the idempotent elements of $*$ in $[0,1)$ (denoted by $k^*$):
\begin{enumerate}[label={\rm(\arabic*)}]
\item If $k^*=1$, then an approach space is $*$-metrizable if and only if it is probabilistic metrizable with respect to the \emph{minimum t-norm} on $[0,1]$.
\item If $k^*<1$, then an approach space is $*$-metrizable if and only if it is probabilistic metrizable with respect to the \emph{product t-norm} on $[0,1]$.
\end{enumerate}
As examples, in Section \ref{Examples} we prove the following:
\begin{itemize}
\item Metric approach spaces are $*$-metrizable for every continuous t-norm $*$ on $[0,1]$ (Proposition \ref{metric-app-*-metrizable}).
\item Approach spaces generated by metrizable topological spaces are $*$-metrizable for every continuous t-norm $*$ on $[0,1]$ (Proposition \ref{topological-app-*-metrizable}).
\item If an approach space is $*$-metrizable for every continuous t-norm $*$ on $[0,1]$, then it is uniform (Corollary \ref{uniform}).
\item If an approach space is $*$-metrizable for some continuous t-norm $*$ on $[0,1]$, then it is locally countable (Proposition \ref{firstcountable}).
\item Approach spaces can be classified by their probabilistic metrizability (Proposition \ref{classification}).
\end{itemize}


\section{Preliminaries}


Throughout this paper, let $[0,\infty]$ denote the extended non-negative real line, with the usual addition ``$+$'' and subtraction ``$-$'' extended via
\[t+\infty=\infty+t=\infty\quad\text{and}\quad\infty-t=\begin{cases}
\infty & \text{if}\ t<\infty,\\
0 & \text{if}\ t=\infty
\end{cases}\]
to $[0,\infty]$. A map 
\[\phi\colon [0,\infty]\to[0,1]\] 
is called a \emph{distance distribution} (see \cite[Definition 4.3.1]{Schweizer1983}) if
\begin{enumerate}[label=(D\arabic*)]
\item $\phi(0)=0$, $\phi(\infty)=1$,
\item $\phi$ is monotone, and
\item $\phi$ is left-continuous on $(0,\infty)$.
\end{enumerate}
In particular, we have
\[\phi(t)=\sup\limits_{s<t}\phi(s)\]
for all $t\in[0,\infty)$.

It is easy to verify the following lemma, which gives a canonical procedure of generating distance distributions:

\begin{lem} \label{monotone-dd}
For every monotone map $\phi\colon[0,\infty]\to[0,1]$, the map
\[\psi\colon[0,\infty]\to[0,1],\quad\psi(t)=\begin{cases}
\sup\limits_{s<t}\phi(s) & \text{if}\ t<\infty,\\
1 & \text{if}\ t=\infty
\end{cases}\]
is a distance distribution. 
\end{lem}

Recall that a binary operation $*$ on an interval $[a,b]$ is a \emph{continuous t-norm}  \cite{Klement2000,Alsina2006} if
\begin{itemize}
\item $([a,b],*,b)$ is a commutative monoid,
\item $p*q\leq p'*q'$ if $p\leq p'$ and $q\leq q'$ in $[a,b]$, and
\item $*\colon[a,b]^2\to[a,b]$ is a continuous function (with respect to the usual topology).
\end{itemize}

It is obvious that $p*q\leq\min\{p,q\}$ for all $p,q\in[a,b]$; in particular, $a*q=a$ for all $q\in[a,b]$. An element $q\in[a,b]$ is \emph{idempotent} if $q*q=q$ and, when it is non-idempotent, it necessarily holds that $q*q<q$. 

\begin{lem} \label{idempotent}
Let $*$ be a continuous t-norm on $[a,b]$. If $q$ is an idempotent element, then $p*q=\min\{p,q\}$ for all $p\in[a,b]$.
\end{lem}

\begin{proof}
If $p\geq q$, then $q=q*q\leq p*q$, and thus $p*q=q$. Suppose that $p\leq q$. Since $a*q=a$ and $q*q=q$, the continuity of $*$ guarantees the existence of $p'\in[a,q]$ such that $p'*q=p$. It follows that
\[p*q=p'*q*q=p'*q=p.\qedhere\]
\end{proof}

Given a continuous t-norm $*$ on $[a,b]$, define
\begin{equation} \label{q-def}
q^-=\sup\{p\in[a,q]\mid p*p=p\}
\end{equation}
for all $q\in[a,b]$. It is clear that $q^-$ is the largest idempotent element in $[a,q]$.

\begin{lem} \label{q-wedge}
$(p*q)^-=\min\{p^-, q^-\}$ for all $p,q\in[a,b]$.
\end{lem}

\begin{proof}
$(p*q)^-\leq\min\{p^-, q^-\}$ follows immediately from $p*q\leq \min\{p,q\}$. Conversely, by Lemma \ref{idempotent},
\[\min\{p^-,q^-\}=p^-*q^-\leq p*q,\]
and consequently $\min\{p^-, q^-\}\leq(p*q)^-$.
\end{proof}

It is well known \cite{Faucett1955,Mostert1957,Klement2000,Klement2004b,Alsina2006} that every continuous t-norm $*$ on the unit interval $[0,1]$ can be written as an \emph{ordinal sum} of three basic t-norms, i.e., the minimum, the product, and the {\L}ukasiewicz t-norm:
\begin{itemize}
\item The \emph{minimum t-norm} $*_M$: $p *_M q=\min\{p,q\}$ for all $p,q\in[0,1]$.
\item The \emph{product t-norm} $*_P$: $p *_P q=pq$ for all $p,q\in[0,1]$.
\item The \emph{{\L}ukasiewicz t-norm} $*_{\L}$: $p *_{\L}q=\max\{0,p+q-1\}$ for all $p,q\in[0,1]$.
\end{itemize}
Explicitly, we say that continuous t-norms $([a_1,b_1],*_1)$ and $([a_2,b_2],*_2)$ are isomorphic if there exists an order-isomorphism 
\[\phi\colon[a_1,b_1]\to[a_2,b_2]\] 
such that
\[\phi\colon([a_1,b_1],*_1,b_1)\to([a_2,b_2],*_2,b_2)\]
is an isomorphism of commutative monoids. Then:

\begin{lem} \label{t-norm-rep} {\rm\cite{Klement2000,Klement2004b,Alsina2006}}
For each continuous t-norm $([0,1],*)$, the set of non-idempotent elements of $*$ in $[0,1]$ is a union of countably many pairwise disjoint open intervals 
\[\{(a_i,b_i)\mid 0<a_i<b_i<1,\ i\in I,\ I\ \text{is countable}\},\]
and for each $i\in I$, the continuous t-norm $([a_i,b_i],*)$ obtained by restricting $*$ to $[a_i,b_i]$ is either isomorphic to the product t-norm $([0,1],*_P)$ or isomorphic to the {\L}ukasiewicz t-norm $([0,1],*_{\L})$.
\end{lem}

The \emph{convolution} of distance distributions $\phi,\psi\colon [0,\infty]\to[0,1]$ with respect to a continuous t-norm $*$ on $[0,1]$ is defined as a map
\[\phi\otimes_*\psi\colon [0,\infty]\to[0,1],\quad(\phi\otimes_*\psi)(t)=\begin{cases}
\sup\limits_{r+s=t}\phi(r)*\psi(s) & \text{if}\ t<\infty,\\
1 & \text{if}\ t=\infty.
\end{cases}\]
It is obvious that the distance distribution
\[\kappa\colon[0,\infty]\to[0,1],\quad\kappa(t)=\begin{cases}
0 & \text{if}\ t=0,\\
1 & \text{if}\ t>0
\end{cases}\]
is the neutral element for the convolution $\otimes_*$. A \emph{probabilistic metric space} (see \cite[Definition 8.1.1]{Schweizer1983}) is a set $X$ equipped with a map 
\[\alpha\colon X\times X\times[0,\infty]\to[0,1],\] 
such that 
\begin{enumerate}[label=(P\arabic*)]
\item \label{PM-def:dd} $\alpha(x,y,-)\colon[0,\infty]\to[0,1]$ is a distance distribution,
\item \label{PM-def:unit} $\alpha(x,x,-)=\kappa$,
\item \label{PM-def:sym} $\alpha(x,y,r)=\alpha(y,x,r)$,
\item \label{PM-def:sep} $\alpha(x,y,-)=\kappa\implies x=y$, and
\item \label{PM-def:tran} $\alpha(y,z,r)*\alpha(x,y,s)\leq\alpha(x,z,r+s)$
\end{enumerate}
for all $x,y,z\in X$, $r,s\in[0,\infty]$, where the value
\[\alpha(x,y,r)\]
is interpreted as the probability that the distance between $x$ and $y$ is less than $r$. A map $f\colon (X,\alpha)\to(Y,\beta)$ between probabilistic metric spaces is \emph{non-expansive} if
\[\alpha(x,x',t)\leq\beta(f(x),f(x'),t)\]
for all $x,x'\in X$, $t\in[0,\infty]$. The category of probabilistic metric spaces (with respect to a continuous t-norm $*$ on $[0,1]$) and non-expansive maps is denoted by
\[\ProbMet_*.\]
Let $\Met$ denote the category of classical metric spaces (in the sense of Fr{\' e}chet \cite{Frechet1906}) and non-expansive maps. There is a canonical embedding functor
\begin{equation} \label{sX}
\chi\colon\Met\to\ProbMet_*
\end{equation}
that sends each metric space $(X,d)$ to the probabilistic metric space $(X,\alpha_d)$ with
\begin{equation} \label{alpha_d-def}
\alpha_d\colon X\times X\times [0,\infty]\to [0,1],\quad \alpha_d(x,y,t)=\begin{cases}
0 & \text{if}\ t\leq d(x,y),\\
1 & \text{if}\ t>d(x,y).
\end{cases}
\end{equation}

\begin{lem} \label{PM-M-P-L}
\begin{enumerate}[label={\rm(\arabic*)}]
\item \label{PM-M-P-L:M} If $(X,\alpha)$ is a probabilistic metric space with respect to $*_M$, then $(X,\alpha)$ is also a probabilistic metric space with respect to every continuous t-norm $*$ on $[0,1]$.
\item \label{PM-M-P-L:P-L} If $(X,\alpha)$ is a probabilistic metric space with respect to $*_P$, then $(X,\alpha)$ is also a probabilistic metric space with respect to $*_{\L}$.
\end{enumerate}
\end{lem}

\begin{proof}
Note that for any $p,q\in[0,1]$, we have $p*q\leq\min\{p,q\}$ and $p*_{\L}q\leq pq$, where the latter follows from $pq-p*_{\L}q=pq-p-q+1=(p-1)(q-1)\geq 0$. Hence, the conclusions follow from \ref{PM-def:tran}.
\end{proof}

Let $\mathbf{2}^X$ denote the powerset of a set $X$. An \emph{approach space} \cite{Lowen1989b,Lowen2015} is a set $X$ equipped with a map 
\[\delta\colon X\times\mathbf{2}^X\to[0,\infty],\] 
such that
\begin{enumerate}[label=(A\arabic*)]
\item \label{app-x} $\delta(x,\{x\})=0$,
\item \label{app-n} $\delta(x,\varnothing)=\infty$,
\item \label{app-u} $\delta(x,A\cup B)=\min\{\delta(x,A),\delta(x,B)\}$, and
\item \label{app-t} $\delta(x,A)\leq\sup\limits_{b\in B}\delta(b,A)+\delta(x,B)$
\end{enumerate}
for all $x\in X$, $A,B\in\mathbf{2}^X$. A map $f\colon (X,\delta)\to(Y,\theta)$ between approach spaces is a \emph{contraction} if 
\[\delta(x,A)\geq\theta(f(x),f(A))\]
for all $x\in X$, $A\in\mathbf{2}^X$. We denote by
\[\App\]
the category of approach spaces and contractions. It is well known (see, e.g., \cite[Theorem 2.2.4]{Lowen2015}) that the category $\Top$ of topological spaces and continuous maps is a coreflective subcategory of $\App$, with the coreflection
\begin{equation} \label{sC}
\sC\colon\App\to\Top
\end{equation}
sending each approach space $(X,\delta)$ to $(X,\cl_{\delta})$, where 
\begin{equation} \label{cl_delta-def}
\cl_{\delta}\colon\mathbf{2}^X\to\mathbf{2}^X,\quad\cl_{\delta}(A)=\{x\in X\mid\delta(x,A)=0\}
\end{equation}
is the topological closure operator induced by the approach structure $\delta$.

\section{Approach structures on probabilistic metric spaces}

Let $(X,\alpha)$ be a probabilistic metric space with respect to a continuous t-norm $*$ on $[0,1]$. Define a map
\begin{equation} \label{delta_alpha-def}
\delta_{\alpha}\colon X\times\mathbf{2}^X\to[0,\infty],\quad\delta_{\alpha}(x,A)=\inf\{r\in[0,\infty]\mid\sup\limits_{a\in A}\alpha(x,a,r)=1\}.
\end{equation}

\begin{prop} \label{alpha-delta-approach}
$(X,\delta_{\alpha})$ is an approach space.
\end{prop}

\begin{proof}
$(X,\delta_{\alpha})$ clearly satisfies \ref{app-x} and \ref{app-n}. For \ref{app-u} and \ref{app-t}, let $x\in X$, $A,B\in\mathbf{2}^X$. First,
\begin{align*}
    \delta_{\alpha}(x,A\cup B) & = \inf\{r\in[0,\infty]\mid\sup\limits_{y\in A\cup B}\alpha(x,y,r)=1\}\\
    & = \inf\{r\in[0,\infty]\mid\sup\limits_{a\in A}\alpha(x,a,r)=1\ \text{or}\ \sup\limits_{b\in B}\alpha(x,b,r)=1\}\\
    & = \min\{\inf\{r\in[0,\infty]\mid\sup\limits_{a\in A}\alpha(x,a,r)=1\},\inf\{r\in[0,\infty]\mid\sup\limits_{b\in B}\alpha(x,b,r)=1\}\}\\
    & = \min\{\delta_{\alpha}(x,A),\delta_{\alpha}(x,B)\}.
\end{align*}
Second, for each $t\in(0,\infty)$, we show that $\delta_{\alpha}(x,A)\leq t$ whenever $\sup\limits_{b\in B}\delta_{\alpha}(b,A)+\delta_{\alpha}(x,B)<t$. In this case, there exist $r,s\in[0,\infty)$ such that
\begin{itemize}
\item $r+s=t$,
\item $\delta_{\alpha}(b,A)<r$ for all $b\in B$, and
\item $\delta_{\alpha}(x,B)<s$.
\end{itemize}
For any $p\in(0,1)$, let $q\in(p,1)$. Then the continuity of $*$ guarantees the existence of $k\in(0,1)$ such that $p<q*k$. Note that $\delta_{\alpha}(x,B)<s$ and the monotonicity of $\alpha(x,b,-)$ guarantees that $\sup\limits_{b\in B}\alpha(x,b,s)=1$. Thus there exists $b\in B$ such that $q<\alpha(x,b,s)$. Similarly, from $\sup\limits_{a\in A}\alpha(b,a,r)=1$ we may find $a\in A$ such that $k<\alpha(b,a,r)$. It follows that 
\[\alpha(x,a,t)\geq\alpha(b,a,r)*\alpha(x,b,s)\geq q*k>p,\]
and consequently $\sup\limits_{a\in A}\alpha(x,a,t)=1$ by the arbitrariness of $p$. Hence $\delta_{\alpha}(x,A)\leq t$.
\end{proof}

\begin{prop}\label{morphism}
Let $(X,\alpha),(Y,\beta)$ be probabilistic metric spaces with respect to a continuous t-norm $*$ on $[0,1]$. If $f\colon (X,\alpha)\to (Y,\beta)$ is non-expansive, then $f\colon(X,\delta_{\alpha})\to (Y,\delta_{\beta})$ is a contraction. 
\end{prop}

\begin{proof}
Let $x\in X$ and $A\subseteq X$. For each $r\in[0,\infty]$, since $f\colon (X,\alpha)\to (Y,\beta)$ is non-expansive, we have
\[\alpha(x,y,r)\leq\beta(f(x),f(y),r)\]
for each $y\in X$. Therefore, 
\[\delta_{\alpha}(x,A)=\inf\{r\in[0,\infty]\mid\sup\limits_{a\in A}\alpha(x,a,r)=1\}\geq \inf\{r\in[0,\infty]\mid\sup\limits_{a\in A}\beta(f(x),f(a),r)=1\}=\delta_{\beta}(f(x),f(A)),\]
which implies that $f\colon(X,\delta_{\alpha})\to (Y,\delta_{\beta})$ is a contraction.
\end{proof}

Thus, we obtain a faithful functor 
\[\Gamma\colon\ProbMet_*\to\App.\]
Note that the diagram 
\[\bfig
\square/^(->`->`<-`->/<1000,500>[\Met`\Top`\ProbMet_*`\App;`\chi`\sC`\Gamma]
\efig\]
is commutative (see \eqref{sX} and \eqref{sC} for the definitions of $\chi$ and $\sC$). Indeed, given a classical metric space $(X,d)$, it is easy to check that
\[\delta_{\alpha_d}(x,A)=d(x,A):=\inf_{a\in A}d(x,a)\]
for all $x\in X$, $A\subseteq X$; and consequently, the topological closure operator induced by $\delta_{\alpha_d}$ (under the functor $\sC$) is precisely the closure operator induced by the metric $d$. Hence, the approach structure \eqref{delta_alpha-def} on a probabilistic metric space is compatible with the metric topology on a classical metric space.

Moreover, for each probabilistic metric space $(X,\alpha)$, the topological space $\sC\Gamma(X,\alpha)$ is exactly the \emph{strong topology} (see \cite[Definition 12.1.1]{Schweizer1983}) for $(X,\alpha)$:
 
\begin{prop}
The neighborhood system $\{\CN_x\}_{x\in X}$ of the topological space $\sC\Gamma(X,\alpha)$ is given by
\[\CN_x=\{U_x(t)\mid t\in(0,\infty)\},\]
where
\[U_x(t)=\{y\in X\mid \alpha(x,y,t)>1-t\}.\] 
\end{prop} 

\begin{proof}
By \eqref{cl_delta-def} and \eqref{delta_alpha-def}, the topological space $\sC\Gamma(X,\alpha)$ is given by the closure operator 
\[\cl_{\delta_{\alpha}}\colon\mathbf{2}^X\to\mathbf{2}^X,\]
with
\[x\in\cl_{\delta_{\alpha}}(A)\iff\inf\{r\in[0,\infty]\mid\sup\limits_{a\in A}\alpha(x,a,r)=1\}=0\]
for all $A\subseteq X$. Thus
\begin{align*}
x\in\cl_{\delta_\alpha}(A) & \iff\forall r\in(0,\infty)\colon\sup_{a\in A}\alpha(x,a,r)=1\\
& \iff \forall r\in(0,\infty),\ \forall t\in (0,r],\ \exists a\in A\colon\alpha(x,a,r)>1-t\\
& \iff \forall t\in(0,\infty),\ \forall r\in[t,\infty),\ \exists a\in A\colon\alpha(x,a,r)>1-t\\
& \iff \forall t\in(0,\infty),\ \exists a\in A\colon\alpha(x,a,t)>1-t & (\alpha(x,y,-)\ \text{is monotone})\\
& \iff \forall t\in(0,\infty)\colon A\cap U_x(t)\neq\varnothing
\end{align*}
for all $A\subseteq X$. Hence, the closure operator $\cl_{\delta_{\alpha}}$ and the neighborhood system $\{\CN_x\}_{x\in X}$ generates the same topology.
\end{proof}

\section{Probabilistic metrizability of approach spaces}

Let $*$ be a continuous t-norm on $[0,1]$. We say that an approach space $(X,\delta)$ is \emph{probabilistic metrizable with respect to $*$}, or \emph{$*$-metrizable} for short, if there exists a probabilistic metric $\alpha$ with respect to $*$ on $X$, such that
\[\delta=\delta_{\alpha}.\] 


\begin{lem} \label{minimpliesall}
Let $*$ be a continuous t-norm on $[0,1]$. Then every $*_M$-metrizable approach space is $*$-metrizable.
\end{lem}

\begin{proof}
This is an immediate consequence of Lemma \ref{PM-M-P-L}\ref{PM-M-P-L:M}.
\end{proof}

\begin{lem}\label{p=l}
An approach space $(X,\delta)$ is $*_P$-metrizable if and only if it is $*_{\text{\L}}$-metrizable.
\end{lem}

\begin{proof}
The necessity follows immediately from Lemma \ref{PM-M-P-L}\ref{PM-M-P-L:P-L}. For the sufficiency, let $\alpha\colon X\times X\times[0,\infty]\to[0,1]$ be a probabilistic metric with respect to $*_{\text{\L}}$ such that $\delta=\delta_{\alpha}$. Define a map
\[\beta\colon X\times X\times [0,\infty]\to [0,1],\quad\beta(x,y,t)=\begin{cases}
0 & \text{if}\ t=0, \\
e^{\alpha(x,y,t)-1} & \text{if}\ t>0.
\end{cases}\]

First, $\beta(x,y,-)$ is a distance distribution. Indeed, it is clear that $\beta(x,y,0)=0$, $\beta(x,y,\infty)=1$ and $\beta(x,y,-)$ is monotone. Since $\alpha(x,y,-)$ is left-continuous on $(0,\infty)$ and the exponential function is continuous, $\beta(x,y,-)$ is also left-continuous on $(0,\infty)$.

Second, $\beta$ is a probabilistic metric with respect to $*_P$, since $\beta$ clearly satisfies \ref{PM-def:unit}, \ref{PM-def:sym} and \ref{PM-def:sep}, while \ref{PM-def:tran} follows from
\[\beta(y,z,r)*_P\beta(x,y,s)=e^{\alpha(y,z,r)+\alpha(x,y,s)-2}\leq e^{\alpha(y,z,r)*_{\L}\alpha(x,y,s)-1}\leq e^{\alpha(x,z,r+s)-1}=\beta(x,z,r+s)\]
for all $x,y,z\in X$, $r,s\in(0,\infty]$. 

Finally, $\delta_{\alpha} =\delta_{\beta}$, because it is easy to see that
\[\sup\limits_{a\in A}\alpha(x,a,r)=1\iff\sup\limits_{a\in A}\beta(x,a,r)=1\]
for all $x\in X$, $A\subseteq X$, $r\in[0,\infty]$.
\end{proof}

Let $*$ be a continuous t-norm on $[0,1]$. Define 
\[k^*=\sup\{q\in[0,1)\mid q*q=q\};\]
that is, the supremum of the idempotent elements of $*$ in $[0,1)$.

\begin{lem} \label{wedgesup}
Let $*$ be a continuous t-norm on $[0,1]$ with $k^*=1$. Then
\[\sup A=1\iff\sup\{a^-\mid a\in A\}=1\]
for all $A\subseteq [0,1]$.
\end{lem}

\begin{proof}
It suffices to prove the ``only if'' part. Suppose that $\sup A=1$. For any $r<1$, from $k^*=1$ we may find an idempotent element $q$ with $r<q<1$, and from $\sup A=1$ we may find $a\in A$ with $q<a\leq 1$. Thus $r<q\leq a^-\leq 1$. Since $r$ is arbitrary, we deduce that $\sup\{a^-\mid a\in A\}=1$.
\end{proof}

\begin{lem}\label{largesummand}
Let $*_1$ and $*_2$ be continuous t-norms on $[0,1]$. If there exists an idempotent element $q$ of $*_1$ such that the continuous t-norms $([q,1],*_1)$ and $([0,1],*_2)$ are isomorphic, then an approach space $(X,\delta)$ is $*_1$-metrizable if and only if it is $*_2$-metrizable.
\end{lem}

\begin{proof}
Let $\phi\colon([q,1],*_1)\to ([0,1],*_2)$ be an isomorphism of continuous t-norms. For the ``only if'' part, suppose that there is a probabilistic metric $\alpha\colon X\times X\times [0,\infty]\to [0,1]$ with respect to $*_1$ such that $\delta=\delta_{\alpha}$. Define a map
\[\beta\colon X\times X\times [0,\infty]\to [0,1],\quad\beta(x,y,t)=
\begin{cases}
0 & \text{if}\ \alpha(x,y,t)<q,\\
\phi(\alpha(x,y,t)) & \text{if}\ \alpha(x,y,t)\geq q.
\end{cases}\]
We show that $\beta$ is a probabilistic metric with respect to $*_2$ and $\delta_{\beta}=\delta$. 

First, $\beta(x,y,-)$ is a distance distribution. It is clear that $\beta(x,y,0)=0$, $\beta(x,y,\infty)=1$ and $\beta(x,y,-)$ is monotone. For the left-continuity of $\beta(x,y,-)$, let $t\in[0,\infty)$. 
\begin{itemize}
\item If $\alpha(x,y,t)<q$, then $\beta(x,y,t)=0=\sup\limits_{s<t}\beta(x,y,s)$.
\item If $\alpha(x,y,t)=q$, then $\beta(x,y,t)=\phi(q)=0=\sup\limits_{s<t}\beta(x,y,s)$.
\item If $\alpha(x,y,t)>q$, then the left-continuity of $\alpha(x,y,-)$ guarantees that there exists $r<t$ such that $\alpha(x,y,s)>q$ for all $s\in(r,t)$. It follows that
\[\beta(x,y,t)=\phi(\alpha(x,y,t))=\phi(\sup_{s\in(r,t)}\alpha(x,y,s))=\sup_{s\in(r,t)}\phi(\alpha(x,y,s))=\sup_{s<t}\beta(x,y,s).\]
\end{itemize}

Second, $\beta$ is a probabilistic metric with respect to $*_2$. It is clear that $\beta$ satisfies \ref{PM-def:unit}, \ref{PM-def:sym} and \ref{PM-def:sep}. For \ref{PM-def:tran}, let $x,y,z\in X$ and $r,s\in [0,\infty]$.
\begin{itemize}
\item If either $\alpha(y,z,r)<q$ or $\alpha(x,y,s)<q$, then 
\[\beta(y,z,r)*_2\beta(x,y,s)=0\leq\beta(x,z,r+s).\]
\item If $\alpha(y,z,r)\geq q$ and $\alpha(x,y,s)\geq q$, then $\alpha(x,z,r+s)\geq q*_1 q=q$, and consequently
\begin{align*}
\beta(y,z,r)*_2\beta(x,y,s)&=\phi(\alpha(y,z,r))*_2\phi(\alpha(x,y,s))\\
&=\phi(\alpha(y,z,r)*_1\alpha(x,y,s))\\
&\leq\phi(\alpha(x,z,r+s))\\
&=\beta(x,z,r+s).
\end{align*}
\end{itemize}

Finally, $\delta_{\beta}=\delta$. Since $\phi$ is an order-isomorphism, we have 
\[\sup_{a\in A}\beta(x,a,r)=1 \iff \sup_{a\in A}\alpha(x,a,r)=1\]
for all $x\in X$, $A\subseteq X$, $r\in [0,\infty]$. The conclusion thus follows.

Conversely, for the ``if'' part, suppose that there is a probabilistic metric $\beta\colon X\times X\times [0,\infty]\to [0,1]$ with respect to $*_2$ on $X$ such that $\delta=\delta_{\beta}$. Define a map
\[\alpha\colon X\times X\times [0,\infty]\to [0,1],\quad
\alpha(x,y,t)=\begin{cases}
0 & \text{if}\ t=0,\\
\phi^{-1}(\beta(x,y,t)) & \text{if}\ t>0.
\end{cases}\]
We show that $\alpha$ is a probabilistic metric with respect to $*_1$ and $\delta_{\alpha}=\delta$. 

First, $\alpha(x,y,-)$ is a distance distribution, since $\alpha(x,y,0)=0$, $\alpha(x,y,\infty)=1$, and 
\[\alpha(x,y,t)=\phi^{-1}(\beta(x,y,t))=\phi^{-1}(\sup_{s<t}\beta(x,y,s))=\sup_{s<t}\phi^{-1}(\beta(x,y,s))=\sup_{s<t}\alpha(x,y,s)\]
for all $t\in(0,\infty)$.

Second, $\alpha$ is a probabilistic metric with respect to $*_1$, since $\alpha$ clearly satisfies \ref{PM-def:unit}, \ref{PM-def:sym} and \ref{PM-def:sep}, while \ref{PM-def:tran} follows from
\[\alpha(y,z,r)*_1\alpha(x,y,s)=\phi^{-1}(\beta(y,z,r))*_1\phi^{-1}(\beta(x,y,s))=\phi^{-1}(\beta(y,z,r)*_2\beta(x,y,s))\leq\alpha(x,z,r+s),\]
for all $x,y,z\in X$, $r,s\in(0,\infty]$. 

Finally, $\delta_{\alpha}=\delta$ follows immediately from the fact that 
\[\sup_{a\in A}\alpha(x,a,r)=1 \iff \sup_{a\in A}\beta(x,a,r)=1\]
for all $x\in X$, $A\subseteq X$, $r\in [0,\infty]$, which completes the proof.
\end{proof}

\begin{thm} \label{main}
Let $*$ be a continuous t-norm on $[0,1]$.
\begin{enumerate}[label={\rm(\arabic*)}]
\item \label{main:=1} If $k^*=1$, then an approach space $(X,\delta)$ is $*$-metrizable if and only if it is $*_M$-metrizable.
\item \label{main:<1} If $k^*<1$, then an approach space $(X,\delta)$ is $*$-metrizable if and only if it is $*_P$-metrizable.
\end{enumerate}
\end{thm}

\begin{proof}
\ref{main:=1} The sufficiency is an immediate consequence of Lemma \ref{minimpliesall}. Conversely, suppose that there is a probabilistic metric $\alpha\colon X\times X\times [0,\infty]\to [0,1]$ with respect to $*$ on $X$ such that $\delta=\delta_{\alpha}$. Define a map
\[\beta\colon X\times X\times [0,\infty]\to [0,1],\quad\beta(x,y,t)=\begin{cases}
\sup\limits_{s<t}\alpha(x,y,s)^- & \text{if}\ t<\infty,\\
1 & \text{if}\ t=\infty.
\end{cases}\]
We show that $\beta$ is a probabilistic metric with respect to $*_M$ on $X$ and $\delta=\delta_{\beta}$.

First, the monotonicity of the map $\alpha(x,y,-)^-\colon[0,\infty]\to[0,1]$ and Lemma \ref{monotone-dd} guarantee that $\beta(x,y,-)$ is a distance distribution.

Second, $\beta$ is a probabilistic metric with respect to $*_M$ on $X$, since $\beta$ clearly satisfies \ref{PM-def:unit}, \ref{PM-def:sym} and \ref{PM-def:sep}, while \ref{PM-def:tran} follows from
\begin{align*}
\beta(y,z,r)*_M\beta(x,y,s) & = (\sup\limits_{a<r}\alpha(y,z,a)^-)*_M(\sup\limits_{b<s}\alpha(x,y,b)^-)\\
& = \sup\limits_{a<r,\,b<s}\min\{\alpha(y,z,a)^-,\alpha(x,y,b)^-\}\\
& = \sup\limits_{a<r,\,b<s}(\alpha(y,z,a)*\alpha(x,y,b))^- & \text{(Lemma \ref{q-wedge})}\\
& \leq \sup\limits_{a<r,\,b<s}\alpha(x,y,a+b)^-\\
& \leq \sup\limits_{t<r+s}\alpha(x,y,t)^-\\
& = \beta(x,z,r+s)
\end{align*}
for all $x,y,z\in X$, $r,s\in[0,\infty)$. 

Finally, $\delta=\delta_{\beta}$. Indeed,
\begin{align*}
\delta_{\beta}(x,A) & = \inf\{r\in[0,\infty]\mid \sup_{a\in A}\beta(x,a,r)=1\}\\
& = \inf\{r\in[0,\infty]\mid \sup_{a\in A}\sup_{s<r}\alpha(x,a,s)^-=1\}\\
& = \inf\{r\in[0,\infty]\mid \sup_{a\in A}\sup_{s<r}\alpha(x,a,s)=1\} & \text{(Lemma \ref{wedgesup})}\\
& = \inf\{r\in[0,\infty]\mid \sup_{a\in A}\alpha(x,a,r)=1\}\\
& = \delta(x,A)
\end{align*}
for all $x\in X$, $A\subseteq X$, as desired.

\ref{main:<1} By Lemma \ref{t-norm-rep}, the continuous t-norm $([k^*,1],*)$ is either isomorphic to $([0,1],*_P)$ or $([0,1],*_\text{\L})$. 
\begin{itemize}
\item If $([k^*,1],*)$ is isomorphic to $([0,1],*_P)$, then it follows from Lemma \ref{largesummand} that $(X,\delta)$ is $*$-metrizable if and only if it is $*_P$-metrizable. 
\item If $([k^*,1],*)$ is isomorphic to $([0,1],*_\text{\L})$, then it follows from Lemmas \ref{p=l} and \ref{largesummand} that $(X,\delta)$ is $*$-metrizable if and only if it is $*_P$-metrizable. \qedhere
\end{itemize}
\end{proof}

\section{Examples} \label{Examples}

\subsection{Metric approach spaces}

Given a (classical) metric space $(X,d)$, the map 
\[\delta_d\colon X\times 2^X\lra[0,\infty],\quad
\delta_d(x,A)=
\begin{cases}
\inf\limits_{y\in A}d(x,y) & \text{if}\ A\neq\varnothing,\\
\infty & \text{if}\ A=\varnothing
\end{cases}\] 
defines an approach structure on $X$. Such approach spaces are called \emph{metric approach spaces} \cite{Lowen1997,Lowen2015}.

\begin{prop} \label{metric-app-*-metrizable}
Metric approach spaces are $*$-metrizable for every continuous t-norm $*$ on $[0,1]$.
\end{prop}

\begin{proof}
Let $(X,\delta)$ be an approach space such that $\delta=\delta_d$ for some metric $d$ on $X$. Let $\alpha_d$ be the probabilistic metric on $X$ induced by $d$ (see \eqref{alpha_d-def}). Then the $*$-metrizability of $(X,\delta)$ follows from
\begin{align*}
\delta_{\alpha_d}(x,A) & = \inf\{r\in[0,\infty]\mid\sup_{a\in A}\alpha_d(x,a,r)=1\}\\
& = \inf\{r\in[0,\infty]\mid\exists\ a\in A\colon \alpha_d(x,a,r)=1\}\\
& = \inf\{r\in[0,\infty]\mid\exists\ a\in A\colon r>d(x,a)\}\\
& = \inf_{a\in A}d(x,a)\\
& = \delta(x,A)
\end{align*}
for all $x\in X$, $A\subseteq X$.
\end{proof}

\subsection{Topological approach spaces}

Given a topological space $(X,\CT)$, the map 
\[\delta_{\CT}\colon X\times 2^X\lra [0,\infty],\quad
\delta_{\CT}(x,A)=
\begin{cases}
0 & \text{if}\ x\ \text{is in the closure of}\ A,\\
\infty & \text{otherwise}
\end{cases}\]
defines an approach structure on $X$. Such approach spaces are called \emph{topologically generated}.

\begin{prop} \label{topological-app-*-metrizable}
If $(X,\CT)$ is a metrizable topological space, then $(X,\delta_{\CT})$ is $*$-metrizable for every continuous t-norm $*$ on $[0,1]$.
\end{prop}

\begin{proof}
By Lemma \ref{minimpliesall}, it suffices to show that $(X,\delta_{\CT})$ is $*_M$-metrizable. Let $d$ be a metric on $X$ such that $x$ is in the closure of $A$ if and only if 
\[\inf_{a\in A}d(x,a)=0\]
for all $x\in X$ and $A\subseteq X$. Let $\phi_0=\kappa$, 
and for each $t\in(0,\infty)$, define a map
\[\phi_t\colon [0,\infty]\to[0,1],\quad \phi_t(x)=
\begin{cases}
1-e^{-\frac{x}{t}} & \text{if}\  x\in[0,\infty),\\
1 & \text{if}\  x=\infty.
\end{cases}
\]
Then $\phi_t$ is clearly a distance distribution for all $t\in[0,\infty)$. We claim that the map 
\[\alpha\colon X\times X\times[0,\infty]\to[0,1],\quad \alpha(x,y,r)=\phi_{d(x,y)}(r)\]
is a probabilistic metric with respect to $*_M$ on $X$. Indeed, since $d$ is a metric, $\alpha$ clearly satisfies \ref{PM-def:unit}, \ref{PM-def:sym} and \ref{PM-def:sep}. For \ref{PM-def:tran}, we have to check that
\begin{equation} \label{topological-app-*-metrizable:M-tran}
\alpha(y,z,r)*_M\alpha(x,y,s)\leq\alpha(x,z,r+s)
\end{equation}
for all $x,y,z\in X$ and $r,s\in[0,\infty)$. Note that \eqref{topological-app-*-metrizable:M-tran} clearly holds when $d(y,z)=0$ or $d(x,y)=0$. Suppose that $d(y,z)\neq 0$ and $d(x,y)\neq 0$. Since 
\[\min\Big\{\frac{r}{d(y,z)},\frac{s}{d(x,y)}\Big\}\leq\frac{r+s}{d(y,z)+d(x,y)}\leq\frac{r+s}{d(x,z)},\]
we have 
\[\alpha(y,z,r)*_M\alpha(x,y,s)=\min\{1-e^{-\frac{r}{d(y,z)}},1-e^{-\frac{s}{d(x,y)}}\}\leq1-e^{-\frac{r+s}{d(x,z)}}=\alpha(x,z,r+s),\]
as desired. Finally, we show that 
\[\delta_{\CT}(x,A)=\delta_\alpha(x,A)\]
for all $x\in X$ and $A\subseteq X$. Let $p=\inf\limits_{a\in A}d(x,a)$. There are two cases:
\begin{itemize}
\item If $p>0$, then $x$ is not in the closure of $A$, and thus $\delta_{\CT}(x,A)=\infty$. On the other hand, since
\[\sup_{a\in A}\alpha(x,a,r)=\sup_{a\in A}\phi_{d(x,a)}(r)=\sup_{a\in A}(1-e^{-\frac{r}{d(x,a)}})=1-e^{-\frac{r}{p}}<1\]
for all $r\in(0,\infty)$, we have $\delta_\alpha(x,A)=\inf\{r\in[0,\infty]\mid\sup\limits_{a\in A}\alpha(x,a,r)=1\}=\inf\{\infty\}=\infty$.
\item If $p=0$, then $\delta_{\CT}(x,A)=0$. On the other hand, it is easy to check that
\[\sup_{a\in A}\alpha(x,a,r)=\sup_{a\in A}\phi_{d(x,a)}(r)=1\]
for all $r\in(0,\infty)$. Hence $\delta_\alpha(x,A)=\inf\{r\in[0,\infty]\mid\sup\limits_{a\in A}\alpha(x,a,r)=1\}=0$, which completes the proof. \qedhere
\end{itemize}
\end{proof}

\subsection{Locally countable approach spaces}

Let $(X,\alpha)$ be a probabilistic metric space. For each $x\in X$ and $n\geq 1$, define a map
\[\lambda_{x,n}\colon X\to [0,\infty],\quad\lambda_{x,n}(y)=\inf\Big\{r\in[0,\infty]\mathrel{\Big|} \alpha(x,y,r)> 1-\dfrac{1}{n}\Big\}.\]

\begin{lem} \label{openball}
Let $(X,\alpha)$ be a probabilistic metric space. Then
\[\delta_{\alpha}(x,A)=\sup_{n\geq 1}\inf_{a\in A}\lambda_{x,n}(a)\]
for all $x\in X$, $A\subseteq X$.
\end{lem}

\begin{proof}
For each $n\geq 1$ and $r\in[0,\infty]$, if $\sup\limits_{a\in A}\alpha(x,a,r)=1$, then there exists $a\in A$ such that $\alpha(x,a,r)> 1-\dfrac{1}{n}$. It follows that 
\[\inf_{a\in A}\lambda_{x,n}(a)=\inf_{a\in A}\inf\Big\{r\in[0,\infty]\mathrel{\Big|} \alpha(x,a,r)>1-\dfrac{1}{n}\Big\}\leq\inf\{r\in[0,\infty]\mid \sup_{x\in A}\alpha(x,a,r)=1\}=\delta_{\alpha}(x,A),\]
and consequently
\[\sup_{n\geq 1}\inf_{a\in A}\lambda_{x,n}(a)\leq \delta_{\alpha}(x,A).\] 
Conversely, we show that for any $s\in(0,\infty)$, $\sup\limits_{n\geq 1}\inf\limits_{a\in A}\lambda_{x,n}(a)<s$ necessarily implies $\delta_{\alpha}(x,A)\leq s$. In this case, for any $n\geq 1$, we have
\[\inf\limits_{a\in A}\lambda_{x,n}(a)=\inf_{a\in A}\inf\Big\{r\in[0,\infty]\mathrel{\Big|}\alpha(x,a,r)> 1-\dfrac{1}{n}\Big\}<s,\]
which guarantees the existence of $a\in A$ such that $\alpha(x,a,s)>1-\dfrac{1}{n}$ (because $\alpha(x,a,-)$ is monotone). It follows that $\sup\limits_{a\in A}\alpha(x,a,s)=1$, and consequently
\[\delta_{\alpha}(x,A)=\inf\{r\in[0,\infty]\mid\sup\limits_{a\in A}\alpha(x,a,r)=1\}\leq s.\qedhere\]
\end{proof}

Let $(X,\delta)$ be an approach space. For each $x\in X$, let
\[\mathcal{A}(x)=\{\phi\colon X\to[0,\infty]\mid \forall A\subseteq X\colon\inf_{a\in A}\phi(a)\leq\delta(x,A)\}.\]
The collection $\{\mathcal{A}(x)\}_{x\in X}$ is called the associated \emph{approach system} of $(X,\delta)$. 
A collection $\{\mathcal{B}(x)\}_{x\in X}$ is a \emph{basis} for the approach system if 
\begin{enumerate}[label=(B\arabic*)]
\item \label{basis:subset} $\mathcal{B}(x)\subseteq\mathcal{A}(x)$ for all $x\in X$,
\item \label{basis:directed} For any $\phi_1,\phi_2\in\mathcal{B}(x)$, there exists $\phi\in\mathcal{B}(x)$ such that $\phi_1,\phi_2\leq\phi$ (under the pointwise order).
\item \label{basis:dominate} Each $\phi\in\mathcal{A}(x)$ is \emph{dominated} by $\mathcal{B}(x)$ in the sense that
\[\forall \epsilon>0,\ \forall\omega<\infty,\ \exists\psi\in\mathcal{B}(x),\ \forall y\in X\colon\min\{\phi(y),\omega\}\leq\psi(y)+\epsilon.\]
\end{enumerate}
$(X,\delta)$ is said to be \emph{locally countable} \cite{Lowen2015} if the associated approach system $\{\mathcal{A}(x)\}_{x\in X}$ has a basis $\{\mathcal{B}(x)\}_{x\in X}$ such that each $\mathcal{B}(x)$ $(x\in X)$ is countable.

\begin{prop} \label{firstcountable}
Every probabilistic metric space with respect to a continuous t-norm $*$ on $[0,1]$ is locally countable.
\end{prop}

\begin{proof}
Let $(X,\alpha)$ be a probabilistic metric space and let $\mathcal{B}(x)=\{\lambda_{x,n}\mid n\geq 1\}$ for each $x\in X$. We show that $\{\mathcal{B}(x)\}_{x\in X}$ is a basis for $\{\mathcal{A}(x)\}_{x\in X}$, the associated approach system of $(X,\delta_{\alpha})$. Indeed, \ref{basis:subset} follows soon from Lemma \ref{openball}, and \ref{basis:directed} holds since for any $n,m\geq 1$, $\lambda_{x,n},\lambda_{x,m}\leq\lambda_{x,\max\{m,n\}}$ is an immediate consequence of the monotonicity of $\alpha(x,y,-)$ $(x,y\in X)$. It remains to prove \ref{basis:dominate}. We proceed by contradiction. Suppose that $\phi\in\mathcal{A}(x)$ is not dominated by $\{\lambda_{x,n}\}_{n\geq 1}$. Then there exist $\epsilon>0$, $\omega<\infty$ such that
\[A_n\coloneqq\{a\in X\mid\min\{\phi(a),\omega\}>\lambda_{x,n}(a)+\epsilon\}\neq\varnothing\]
for all $n\geq 1$. However, since $A_n\subseteq A_m$ whenever $m\leq n$, we have
\begin{align*}
\sup_{n\geq 1}\delta_{\alpha}(x,A_n)+\epsilon
& = \sup_{n\geq 1}\sup_{m\geq 1}\inf_{a\in A_n}\lambda_{x,m}(a)+\epsilon & \text{(Lemma \ref{openball})}\\
& \leq \sup_{n\geq 1}\sup_{m\geq 1}\inf_{a\in A_{\max\{n,m\}}}\lambda_{x,\max\{n,m\}}(a)+\epsilon & (A_{\max\{n,m\}}\subseteq A_n\ \text{and}\ \lambda_{x,m}\leq\lambda_{x,\max\{n,m\}})\\
& = \sup_{n\geq 1}\inf_{a\in A_n}\lambda_{x,n}(a)+\epsilon\\
& = \sup_{n\geq 1}\inf_{a\in A_n}(\lambda_{x,n}(a)+\epsilon)\\
& \leq \sup_{n\geq 1}\inf_{a\in A_n}\min\{\phi(a),\omega\} & (\text{definition of}\ A_n)\\
& = \min\{\sup_{n\geq 1}\inf_{a\in A_n}\phi(a),\omega\}\\
& \leq \min\{\sup_{n\geq 1}\delta_{\alpha}(x,A_n),\omega\}, & (\phi\in\mathcal{A}(x))
\end{align*}
which is a contradiction.
\end{proof}

As it is known from \cite[Proposition 3.4.3]{Lowen2015} that a topologically generated approach space is locally countable if and only if its underlying topology is first-countable, the following corollary follows immediately from Proposition \ref{firstcountable}:

\begin{cor} \label{notfirstcountable}
Let $(X,\delta)$ be a topologically generated approach space. If the underlying topology of $(X,\delta)$ is not first-countable, then $(X,\delta)$ is not $*$-metrizable for any continuous t-norm $*$ on $[0,1]$.
\end{cor}

\subsection{Uniform approach spaces}

A \emph{generalized metric} \cite{Lawvere1973} (also \emph{quasi-metric} \cite{Lowen2015}) on a set $X$ is a map $d\colon X\times X\to [0,\infty]$ satisfying 
\[d(x,x)=0\quad\text{and}\quad d(y,z)+d(x,y)\geq d(x,z)\]
for all $x,y,z\in X$, and it is called \emph{symmetric} if $d(x,y)=d(y,x)$ for all $x,y\in X$. We denote by
\[\CM(X)\]
the family of all generalized metrics on $X$. Let $\rho\in\CM(X)$ and $\CD\subseteq\CM(X)$. We say that $\rho$ is \emph{locally dominated by} $\CD$ if 
\[\forall x\in X,\ \forall\epsilon>0,\ \forall\omega<\infty,\ \exists\ d\in\CD,\ \forall y\in X\colon\min\{\rho(x,y),\omega\}\leq d(x,y)+\epsilon,\]
and
\[\hat{\CD}=\{\rho\in\CM (X)\mid \rho\ \text{is locally dominated by}\ \CD\}\]
is called the \emph{local saturation} of $\CD$.

Let $(X,\delta)$ be an approach space, and let $\{\mathcal{A}(x)\}_{x\in X}$ be the associated approach system of $(X,\delta)$. The collection
\[\CG=\{d\in\CM(X)\mid\forall x\in X\colon d(x,-)\in\CA(x)\}\]
is called the associated \emph{gauge} of $(X,\delta)$. A collection $\CB\subseteq\CG$ is called a \emph{basis} for the associated gauge of $(X,\delta)$ if 
\[\CG=\hat{\CB}.\]
It is known that every generalized metric on $X$ which is locally dominated by $\CG$ belongs to $\CG$ (see \cite[Theorem 1.2.4]{Lowen2015}). An approach space is said to be \emph{uniform} \cite{Lowen2015} if it has a gauge basis consisting of symmetric generalized metrics.

\begin{prop} \label{uniformapproachspace}
Every probabilistic metric space with respect to $*_M$ is a uniform approach space.
\end{prop}

\begin{proof}
Let $(X,\alpha)$ be a probabilistic metric space with respect to $*_M$. For each $n\geq 1$, define a map 
\[d_n\colon X\times X\to[0,\infty],\quad d_n(x,y)=\lambda_{x,n}(y)=\inf\Big\{r\in[0,\infty]\mathrel{\Big|} \alpha(x,y,r)> 1-\dfrac{1}{n}\Big\}.\] 
It is clear that $d_n(x,x)=0$ and $d_n(x,y)=d_n(y,x)$ for all $x,y\in X$. Let $x,y,z\in X$. 
If $r,s\in[0,\infty]$ satisfy 
\[\alpha(y,z,s)>1-\dfrac{1}{n}\quad\text{and}\quad \alpha(x,y,r)>1-\dfrac{1}{n},\] 
then 
\[\alpha(x,z,r+s)\geq\alpha(y,z,s)*_M\alpha(x,y,r)>1-\dfrac{1}{n},\]
and consequently $d_n(y,z)+d_n(x,y)\geq d_n(x,z)$. Thus $d_n$ is a symmetric generalized metric on $X$. We show that 
\[\CB=\{d_n\mid n\geq 1\}\] 
is a basis for the associated gauge $\CG$ of $(X,\delta_\alpha)$. On one hand, $\hat{\CB}\subseteq\CG$ because it is clear that $d_n\in\CG$ for all $n\geq 1$. On the other hand, let $\rho\in\CG$. For any $x\in X$, $\epsilon>0$, $\omega<\infty$, since $\rho(x,-)$ is dominated by $\mathcal{B}(x)=\{\lambda_{x,n}\mid n\geq 1\}$, there exists $n\geq 1$ such that
\[\min\{\rho(x,y),\omega\}\leq \lambda_{x,n}(y)+\epsilon=d_n(x,y)+\epsilon,\]
which means that $\rho$ is locally dominated by $\CB$, and consequently $\rho\in\hat{\CB}$. Hence $\CG\subseteq\hat{\CB}$.
\end{proof}

\begin{cor} \label{uniform}
If an approach space is $*$-metrizable for every continuous t-norm $*$ on $[0,1]$, then it is uniform.
\end{cor}

\subsection{Classifications of approach spaces by probabilistic metrizability}

The following two lemmas are easily verified with the aid of Theorem \ref{main}:

\begin{lem} \label{k*=1}
If an approach space $(X,\delta)$ is not $*_0$-metrizable for some continuous t-norm $*_0$ on $[0,1]$, then it is not $*_M$-metrizable, and consequently not $*$-metrizable for any continuous t-norm $*$ on $[0,1]$ with $k^*=1$.
\end{lem}

\begin{proof}
Since $(X,\delta)$ is not $*_0$-metrizable, it follows from Lemma \ref{minimpliesall} that $(X,\delta)$ cannot be $*_M$-metrizable. Hence, by Theorem \ref{main}\ref{main:=1}, $(X,\delta)$ is not $*$-metrizable for any continuous t-norm $*$ on $[0,1]$ with $k^*=1$.
\end{proof}

\begin{lem} \label{k*<1}
If an approach space $(X,\delta)$ is $*_0$-metrizable for some continuous t-norm $*_0$ on $[0,1]$, then it is $*_P$-metrizable, and consequently $*$-metrizable for every continuous t-norm $*$ on $[0,1]$ with $k^*<1$.
\end{lem}

\begin{proof}
If $k^{*_0}=1$, then $(X,\delta)$ is $*_M$-metrizable by Theorem \ref{main}\ref{main:=1}, and thus $*$-metrizable for every continuous t-norm $*$ on $[0,1]$ by Lemma \ref{minimpliesall}.

If $k^{*_0}<1$, then from Theorem \ref{main}\ref{main:<1} we deduce that $(X,\delta)$ is $*_P$-metrizable, and consequently $*$-metrizable for every continuous t-norm $*$ on $[0,1]$ with $k^*<1$.
\end{proof}

Therefore, we obtain the following classifications of approach spaces:

\begin{prop} \label{classification}
Let $(X,\delta)$ be an approach space. Then exactly one of the following three statements is true:
\begin{enumerate}[label={\rm(\arabic*)}]
\item \label{classification:is} $(X,\delta)$ is $*$-metrizable for every continuous t-norm $*$ on $[0,1]$.
\item \label{classification:m} $(X,\delta)$ is $*_P$-metrizable, but not $*_M$-metrizable; in other words, $(X,\delta)$ is $*$-metrizable for every continuous t-norm $*$ on $[0,1]$ with $k^*<1$, but not $*$-metrizable for any continuous t-norm $*$ on $[0,1]$ with $k^*=1$.
\item \label{classification:not} $(X,\delta)$ is not $*$-metrizable for any continuous t-norm $*$ on $[0,1]$.
\end{enumerate}
\end{prop}


We have known that every metric approach space satisfies Proposition \ref{classification}\ref{classification:is} (see Proposition \ref{metric-app-*-metrizable}), and every topologically generated approach space whose underlying topology is not first-countable satisfies Proposition \ref{classification}\ref{classification:not} (see Corollary \ref{notfirstcountable}). However, it seems not easy to construct an approach space that satisfies \ref{classification}\ref{classification:m}. So, we end this paper with the following:

\begin{ques}
Can we find an example of an approach space that is $*_P$-metrizable but not $*_M$-metrizable?
\end{ques}

\section*{Acknowledgement}

The first and the second named authors acknowledge the support of National Natural Science Foundation of China (No. 12071319 and No. 12171342) and the Fundamental Research Funds for the Central Universities (No. 2021SCUNL202). 





\bibliographystyle{abbrv}
\bibliography{lili}

\end{document}